\documentclass[12pt]{article}
\usepackage{amssymb}
\newtheorem{prop}{Proposition}[section]

\newcommand{\subsect}[1]{\subsection {#1}\mbox{}\par}
\newcommand{\qed}{\hfill $\square$\vskip .2cm}

\newcommand{\sect}[1]{\section{#1}\setcounter{equation}{0}}

\def\<{\langle}
\def\>{\rangle}
\def\Epsilon{{\cal E}}
\def\f{\mathfrak f}
\def\g{\mathfrak g}

\def\L{\mathfrak L}
\def\R{\mathbb R}
\def\C{{\cal C}}
\def\A{\mathfrak A}
\def\S{\mathfrak S}
\def\Z{\mathbb Z}

\def\R{\mathbb R}
\def\G{{\cal G}}
\def\b{\bar{b}}
\begin{document}

 \title{Superdiffusivity of occupation-time variance in
 $2$-dimensional asymmetric exclusion processes with density $\rho =
 1/2$}
 \author{
\begin{tabular}[l]{ccc}
Sunder Sethuraman\\
Iowa State University
\end{tabular}
}

 \thispagestyle{empty}
 \maketitle
 \abstract{We compute that the growth of the origin
occupation-time variance up to time $t$ in dimension $d=2$ with respect to asymmetric
simple exclusion in equilibrium with density $\rho = 1/2$ is in a
 certain sense at least
 $t\log (\log t)$ for general rates, and at least
 $t(\log t)^{1/2}$ for rates which are asymmetric only in the
 direction of one of the axes.  These estimates are consistent
 with an important conjecture with respect to the transition function
 and variance of ``second-class'' particles.
 }
 
 \vskip .2cm
 {\sl Abbreviated title}:  Occupation-time variance
in $2$D asymmetric exclusion process with density $1/2$ \\[.15cm]
 {\sl AMS (2000) subject classifications}: Primary 60K35; secondary 60F05.

\sect{Introduction and Results}
It is known that the occupation-time variance at the origin up to
time $t$ in asymmetric exclusion processes in equilibrium is proportional to $t$
times the expected time a second-class particle,
beginning at the origin, spends at the origin, that is, $t\int_0^t
(1-s/t)p_s(0,0)ds$ where $p_s(0,j)$ is the second-class particle
transition function.  Let us now fix the equilibrium density $\rho =
1/2$ so that the mean of the second-class particle at time $t$ vanishes.
Recently, it has been
argued, as the variance of a second-class particle at time $t$ starting initially
at the origin--$\sum j^2 p_t(0,j)$ in this case--is conjectured to be on the order $t^{4/3}$ in $d=1$
\cite{BKS} and proved (for a closely related resolvent quantity) to be at least $t^{5/4}$ \cite{LQSY}, that the
the transition function of the second-class particle decays on
order $t^{-2/3}$ in $d=1$ (cf. equation (4.8) \cite{PS}).  In $d=2$, the second-class
particle variance is conjectured as $O(t(\log t)^{2/3})$ \cite{BKS}
with a proof (for a resolvent quantity when the process rates are
asymmetric only in the direction of one of the axes)
\cite{Yau}.
Perhaps by the same sort
of reasoning as in \cite{PS}, one may claim the second-class
transition function decays as
$t^{-1}(\log t)^{-1/3}$ in $d=2$ (cf. equation (12) \cite{BKS}).
Then, the
occupation-time variance orders should match second-class particle
variance orders
in both $d=1$ and $2$.  We mention also these variance orders
have connections to fluctuation orders of the current across the origin
on which there has been much study (cf. \cite{F}, \cite{PS}, \cite{Yau}).

In this note, we show that the occupation-time variance at the origin
diverges in a sense in $d=2$ when density $\rho = 1/2$ at least as
$t\log (\log  t)$ for general asymmetric rates, and at least as $t(\log
t)^{1/2}$ when the asymmetry is only in the direction of one of the
axes 
(Proposition
\ref{lb1}) and so is consistent with the above discussion.
The methods are 
to link occupation-time variances and certain resolvent $H_{-1}$
norms, and then to use some ``free-particle'' comparisons of H.T. Yau
in the style of Bernardin \cite{B}.  
\vskip .2cm

{\it Model.}  Informally, the simple exclusion process on $Z^d$ is a collection of random
walks which move with jump rates $p(i,j) = p(j-i)$ independently
except in that jumps to occupied vertices are suppressed.  In this
article, we will assume that $p$ is finite-range such that its
symmetrization $(p(\cdot) + p(-\cdot))/2$ is irreducible.  More
formally, let
$\Sigma = \{0,1\}^{Z^d}$ be the configuration space where
a configuration
$\eta=\{\eta_i: i\in Z^d\}$ is a collection of ``occupation''
coordinates where
$\eta_i = 1$ if $i$ is occupied and $\eta_i=0$ otherwise.
The exclusion process is a Markov process $\eta(t)$ evolving on
$\Sigma$ with generator
$$Lf(\eta) \ = \ \Sigma_{i,j}
p(j-i)\eta_i(1-\eta_j)(f(\eta^{ij})-f(\eta)).$$
Here, $\eta^{ij}$ is the configuration obtained from $\eta$ by
interchanging the values at $i$ and $j$.
Let also $T_t$ denote the associated semi-group.  See \cite{Liggett}
for more details.

It is well-known that
there is a family of 
invariant measures $\{P_\rho: 0\leq \rho\leq 1\}$ each of which
concentrate on configurations of a fixed density $\rho$.  These
measures take form as Bernoulli product measures, that is,
$P_\rho$ independently places a particle at each vertex with
probability $\rho$.  Let $E_\rho$ denote expectation with respect to
$P_\rho$.  Denote also by $\<\cdot,\cdot\>_\rho$ and $\|\cdot\|_0$ the
innerproduct and norm on $L^2(P_\rho)$.

We also note one can compute, with respect to $P_\rho$, that the
adjoint $L^*$ is the generator of simple exclusion with reversed jump
rates $p(-\cdot)$.  

\vskip .2cm

{\it Problem and Connection to Second-Class Particles.}
Consider the centered occupation time, say, at the origin up to time $t$,
$A_\rho(t)= \int_0^t (\eta_0(s)-\rho)ds$. 
The problem is to compute the variance of
$A_\rho(t)$ under the equilibrium $P_\rho$.   Let $\sigma^2_t =
E_\rho[ A^2_\rho(t)]$ denote the variance.  We compute, using
stationarity and basic calculations, that
\begin{eqnarray*}
\sigma^2_t &=& 2\int_0^t\int_0^s E_\rho[(\eta_0(s) - \rho)(\eta_0(0) -
\rho)]du ds\\
&=& 2\int_0^t (t-s) E_\rho[(\eta_0(s) - \rho)(\eta_0(u) -
\rho)]ds.
\end{eqnarray*}
To express the kernel further, consider the ``basic coupling''
of two systems, the first starting
under $\xi \sim P_\rho(\cdot|\eta_0=0)$ and the second under $\xi +
\delta_0$, that is with an extra particle at the origin.
Let
$(\xi(t),R(t))\sim \bar{P}$ denote the coupled process where $R(t)$ tracks the
discrepancy or ``second-class'' particle.  The joint generator is
\begin{eqnarray*}
({\bar{L}}f)(\xi,r) &=&   \sum_{i,j\neq r}
p(j-i)\xi_i(1-\xi_j)(f(\xi^{ij},r)-f(\xi,r))\\
&&\ + \sum_i p(-i)\xi_{r-i}(f(\xi^{r-i,r},r-i)-f(\xi,r))\\
&&  \  +\sum_i
p(i)(1-\xi_{r+i})(f(\xi^{r,r+i},r+i)-f(\xi,r)).
\end{eqnarray*}
The first sum refers to jumps not including the discrepancy, while the
second and third sums correspond to jumps of other particles to
the discrepancy position and jumps of the discrepancy itself.  

We have then
\begin{eqnarray*}
&&E_\rho[(\eta_0(s) - \rho)(\eta_0(0) -
\rho)]\\
&&\ \ \ \ \ \ \ =
\rho(1-\rho)P_\rho(\eta_0(s) =1|\eta_0(0)=1) -
P_\rho(\eta_0(s)=1|\eta_0(0)=0)\\
&&\ \ \ \ \ \ \ = \rho(1-\rho)\bar{P}[R(s) = 0]
\end{eqnarray*}
which leads to
$$\lim_{t\rightarrow \infty}\sigma^2_t/t \ = \ \lim_{t\rightarrow
  \infty}
2\rho(1-\rho)\int_0^t (1-s/t)\bar{P}[R(s) = 0]ds \ = \
  2\rho(1-\rho)\int_0^\infty
\bar{P}[R(s) = 0]ds;$$
the notation earlier in the introduction reads now
$p_t(0,j)=\bar{P}(R(t)=j)$.

The
second-class particle process $R(t)$, with respect to its own history, is not Markov except
when the jump rate $p$ is symmetric, in which case, it is a symmetric
random walk.  In general, it is highly dependent on the whole system.
However, one can roughly think of $R(t)$ as some sort of random walk
with mean drift $(1-2\rho)\sum_i ip(i)$.  This drift vanishes exactly
when $p$ is either mean-zero ($\sum_iip(i) = 0$) or $\rho =1/2$, and so one might think
  the process is recurrent exactly in this case
so that
\begin{eqnarray*}
\lim_{t\rightarrow \infty} \sigma^2_t/t &=& \infty \ \ {\rm in \ }
d\leq 2 \ {\rm when \ } p {\rm \ mean-zero \ or \ when \ }\rho = 1/2\\
&<& \infty \ \ {\rm otherwise.}
\end{eqnarray*}
This has been established in dimensions $d\geq 3$, $d=1$, and in $d=2$
when $\rho \neq 1/2$ \cite{K}, \cite{Sclt}, \cite{SS}, \cite{B}.  
Still open it seems is to show the variance is
superdiffusive in $d=2$ when $\rho = 1/2$.

Of key interest is also how fast $\sigma^2_t/t$ diverges in $d\leq
2$ when $p$ mean-zero or
$\rho =1/2$. In fact, it has been shown that $\sigma^2_t \sim t^{3/2}$
and $t\log t$
in $d=1$ and $d=2$ respectively
when $p$ is mean-zero \cite{K}, \cite{Sclt}.  When $p$ has a drift
$(\sum_i i p(i) \neq 0$) and $\rho = 1/2$, as mentioned earlier,
$\sigma^2_t$ is
conjectured to diverge as $t^{4/3}$ and $t(\log
t)^{2/3}$ in $d=1$ and $d=2$ respectively.
Indeed, a lower bound on order $t^{5/4 - }$ has been shown in \cite{B} in
$d=1$.  The main result of this article (Proposition \ref{lb1}) is to compute in $d= 2$ when
$p$ has a general drift and $\rho =1/2$
that 
$$\liminf_{\lambda \rightarrow 0}\frac{1}{\log (|\log
\lambda|)} \int_0^\infty e^{-\lambda t} \bar{P}(R(t) = 0)dt \ > \ 0 $$
or by integrating the second-class transition function twice
$$ \liminf_{\lambda \rightarrow 0}\frac{\lambda^2}{\log(|\log \lambda |)}\int_0^\infty e^{-\lambda t} \sigma^2_t dt
\ > \ 0.$$
When the drift $\sum i p(i)$ is in the direction of one of the axes,
the same result holds with ``$\log|\log \lambda|$'' replaced by
``$|\log \lambda |^{1/2}$.''
Clearly $\sigma^2_t/t$ diverges regardless, 
and moreover a formal Tauberian analogy would 
suggest that $\sigma^2_t$ is at least on order $t\log(\log t)$ in the
general case and $t(\log t)^{1/2}$ in the more special case.

We mention that some rough upper bounds in $d=1,2$ in the
``drift'' case when $\rho = 1/2$ are easy to obtain by a comparison
with the symmetrized process, namely $\sigma^2_t \leq c_1(\rho)t^{3/2}$
in $d=1$ and $\leq c_2(\rho)t\log t$ in $d=2$.  Although well known,
we include them for completeness in Proposition \ref{ub2}.


\vskip .2cm
{\it Variational Formulas.}  The method of proof does not work with second-class
particles, but with certain variational formulas for some resolvent
quantities.  By a local function, we mean a function supported on a
finite number of coordinates.

The generator $L$ can be decomposed into symmetric and anti-symmetric
parts,
$L=S+A$ where $S=(L +L^*)/2$ and $A=(L-L^*)/2$.  Since $L$ is
Markovian, $S$ is a non-positive operator.  Consider the
resolvent operator $(\lambda -L)^{-1}:L^2(P_\rho) \rightarrow
L^2(P_\rho)$ 
well defined for $\lambda>0$--in particular, $(\lambda -L)^{-1} f =
\int_0^\infty e^{-\lambda s}T_s f ds$.  Since the symmetrization of $(\lambda
-L)^{-1}$ has inverse $(\lambda - L^*)(\lambda -S)^{-1}(\lambda -L) =
(\lambda -S) +A^*(\lambda-S)^{-1}A$, we have
the variational formula for $f$ local,
$$\<f,(\lambda-L)^{-1} f\>_\rho \ = \ \sup_{\phi \ {\rm local}} \big\{
2\<f,\phi\>_\rho - \<\phi, (\lambda -S)\phi\>_\rho -\<A\phi, (\lambda
- L)^{-1}A\phi\>_\rho\big \}.$$
Now, as $A^* = -A$ and $A^*(\lambda -S)^{-1}A$ is a non-positive
operator,
we have the easy bound that $\<f,(\lambda-L)^{-1}f\>_\rho$ is bounded
by its ``symmetrization,''
\begin{eqnarray}
\<f,(\lambda-L)^{-1} f\>_\rho & \leq & \sup_{\phi \ {\rm local}} \big\{
2\<f,\phi\>_\rho - \<\phi, (\lambda -S)\phi\>_\rho \big\}\nonumber\\
&=& \<f,(\lambda -S)^{-1}\>_\rho.
\label{symbound}
\end{eqnarray}

\vskip .2cm
{\it Upper bounds.}\
Well known upperbounds on $\sigma^2_t$ follow from two statements
which we include here for completeness.

\begin{prop} 
\label{ub1} There is a universal constant $C_1$ such that 
\begin{eqnarray*}
\sigma^2_t &\leq& C_1t\<\eta_0-\rho, (t^{-1} - L)^{-1}
\eta_0-\rho\>_\rho\\
&\leq&C_1t\<\eta_0-\rho, (t^{-1} - S)^{-1}
\eta_0-\rho\>_\rho.
\end{eqnarray*}
\end{prop}

{\it Proof.} The first line is well-known (with a proof found for
instance in Lemma 3.9 \cite{Sclt}), and the second is
(\ref{symbound}). \qed

\begin{prop}
\label{ub2} In $d\leq 2$, there exists a constant $C_2 = C_2(d,\rho,p)$ where
for large $t$,
$$\<\eta_0-\rho, (t^{-1} - S)^{-1}
\eta_0-\rho\>_\rho \ \leq \ 
\left\{\begin{array}{rl}C_2\sqrt{t} & \
    {\rm in \ } d=1\\
C_2\log t& {\rm in \ } d=2\end{array}\right.$$
and so by Proposition \ref{ub1}, $\sigma_t^2 \leq C_1C_2 t^{3/2}$ in
$d=1$ and $C_1C_2t\log t$ in $d=2$.
\end{prop}

{\it Proof.} This is proved in \cite{K} as follows:  Write $\<\eta_0-\rho, (t^{-1} - S)^{-1}
\eta_0-\rho\>_\rho = \int_0^\infty e^{-s/t}
E_\rho[(\eta_0(s)-\rho)(\eta_0(0)-\rho)]ds = \rho(1-\rho)\int_0^\infty e^{-s/t} \bar{P}(R_s =
0)ds$.  As in the symmetric case
$\bar{P}(R_s = 0) \sim s^{-1/2}$ in $d=1$ and $s^{-1}$ in $d=2$, the estimates follow. \qed

\vskip .2cm
{\it Lower bounds.}
The lowerbounds are through variational formulas.
The following is the main result of this note and is proved in
subsection 2.1.
Let $e_1$ and $e_2$ denote the standard basis in $\R^2$.
\begin{prop}
\label{lb1}
In $d=2$, when $\sum_i ip(i)\neq 0$ and $\rho=1/2$, there is a
constant $C_3=C_3(\rho,p)$ where for all small $\lambda$,
$$2\rho(1-\rho)\int_0^\infty e^{-\lambda t}\bar{P}(R(t) =0)dt \ = \ \<\eta_0-\rho, (\lambda - L)^{-1}
\eta_0-\rho\>_\rho \ \geq \ C_3\log (|\log \lambda)|);$$
when, more specifically, $\sum_i ip(i)=c e_1$ or $c e_2$  is a non-zero multiple of either
$e_1$ or $e_2$, 
$\<\eta_0-\rho, (\lambda - L)^{-1}
\eta_0-\rho\>_\rho  \geq  C_3|\log  \lambda|^{1/2}$.

\end{prop}




\section{Some Preliminaries}
We first give some tools and definitions before going to the proof of
Proposition \ref{lb1} in subsection 2.1.
\vskip .2cm

{\it Comparison Bound.}  We compare
$\<f,(\lambda -L)^{-1}\>_\rho$ with the formula with respect to a
``nearest-neighbor'' operator $L_0$.
Let
$m_i = e_i \cdot \sum i p(i)$ for $i=1,2$.  As
the drift of $p$ is assumed not to vanish, at least one of the $m_i$'s is
not zero.  Without loss of generality, suppose $m_1\neq 0$.

Let $L_0$ be the exclusion generator corresponding to nearest-neighbor
jump rates $p_0(\cdot)$ where 
$$\begin{array}{l}
p_0(e_1) = |m_1|, p_0(e_2)=|m_2|, \ {\rm and \ } p_0(i)=0 {\rm \
  otherwise,\ when \ }m_2\neq 0 \ \ {\rm and }\\
p_0(e_1) = |m_1|, p_0(\pm e_2)=1/4, {\rm \ and \ } p_0(i) = 0 \ {\rm otherwise,
  \ when \ } m_2= 0.
\end{array}
$$

The following is proved in Theorem 2.1 \cite{Scomp}.
\begin{prop}
\label{comp}
There is a constant $C_4=C_4(p)$ where
$$C_4^{-1}\<f,(\lambda -L_0)^{-1}\>_\rho \ \leq \ \<f,(\lambda
-L)^{-1}\>_\rho \ \leq \ C_4\<f,(\lambda -L_0)^{-1}\>_\rho. $$
\end{prop}

\vskip .2cm

{\it Duality.}
Let $\Epsilon$ denote the collection of finite subsets of $Z^2$, and let
$\Epsilon_n$ denote those subsets of cardinality $n$.  Let also $\Psi_B$ be
the function 
$$\Psi_B(\eta) \ =\ \prod_{x\in B} \frac{\eta_x -
  \rho}{\sqrt{\rho(1-\rho)}}$$
where we take $\Psi_\emptyset = 1$ by convention.  One can check that
  $\{\Psi_B: B\in \Epsilon\}$ is Hilbert basis of $L^2(P_\rho)$.  In
  particular, any function $f\in L^2(P_\rho)$ has decomposition
$$f \ = \ \sum_{n\geq 0} \sum_{B\in \Epsilon_n} \f(B)\Psi_B$$
with coefficient $\f:\Epsilon \rightarrow \R$ which in general depends on $\rho$.

Then, for $f,g\in L^2(P_\rho)$,
$$\<\f,\g\>\ := \ \<f,g\>_\rho \ = \ \sum_{B\in \Epsilon}\f(B)\g(B)$$
and $\|\f\|^2 := \|f\|_0^2 = \<f,f\>_\rho$.
Let $\C_n$ be the subspace generated by local functions of degree $n$,
that is functions whose support sets are members of $\Epsilon_n$.

The operators $L$, $S$ and $A$ have counterparts $\L$, $\S$ and $\A$ which act on
``coefficient'' functions $\f$.  These are given in the expressions
$$Lf \ = \ \sum_{B\in \Epsilon} (\L\f)(B)\Psi_B, \ \ Sf \ = \
\sum_{B\in \Epsilon} (\S\f)(B)\Psi_B, \ \ {\rm and \ } Af \ = \
\sum_{B\in \Epsilon} (\A\f)(B)\Psi_B.$$

Let $s$ and $a$ be the symmetric and anti-symmetric parts of $p$,
$s(i)= (p(i) +p(-i))/2$ and $a(i) = (p(i)-p(-i))/2$.  Also for
$B\subset Z^d$, denote
$$B_{x,y} = \left\{\begin{array}{rl}
B\setminus \{x\} \cup \{y\} & \ {\rm when \ } x\in B, y\not\in B\\
B\setminus \{y\} \cup \{x\} & \ {\rm when \ } x\not\in B, y\in B\\
B&\ {\rm otherwise.}\end{array}\right.
$$

Now, of course, $\L = \S +\A$.  
Moreover, the symmetric part $\S$ can be computed as
$$
(\S\f)(B) \ = \ \frac{1}{2}\sum_{x,y\in \Z^d} s(y-x)[\f(B_{x,y}) -
\f(B)].$$
Also, the anti-symmetric part $\A$ can be decomposed into the sum of
three operators which preserve, increase, and decrease the degree of
the function acted upon: $\A = (1-2\rho)\A_0 +
2\sqrt{\rho(1-\rho)}(\A^+ - \A^-)$.
$$\begin{array}{rl}
(\A_0\f)(B) & = \ \sum_{x\in B \atop y\not\in B}
a(y-x)[\f(B_{x,y})-\f(B)]\\
(\A^+\f)(B) & = \ \sum_{x\in B\atop y\in B} a(y-x) \f(B-\{y\})\\
(\A^-\f)(B) &= \ \sum_{x\not\in B\atop y\not \in B} a(y-x)
\f(B\cup \{x\}).
\end{array}
$$ 

We note from the expression that in fact $\A_0$, $\A^+$ and $\A^-$ take a degree $n$ function,
that is say 
$\f:(\Z^d)^n\rightarrow \R$, into respectively a degree $n$, $n+1$ and
$n-1$ function.
It will be helpful to write $\A$ in terms of its ``degree''
actions,
$$\A \ = \ \sum_{n\geq 0} \bigg ( \A_{n,n-1} + \A_{n,n} +
\A_{n,n+1}\bigg )$$
where $\A_{m,n}$ is the part which takes degree $m$ functions to
degree $n$ functions.
Here, by convention $\A_{0,-1}\f = 0$ is the zero function.

At this point, we observe when $\rho=1/2$ that $\A =
2\sqrt{\rho(1-\rho)}(\A^+ - \A^-)$ as the part which preserves
degree vanishes here.  

\vskip .2cm
{\it $H_{1}$ and $H_{-1}$ Spaces.}
Define, for local functions $f$, the $H_1(P_\rho)$ (semi)-norm by $\|f\|^2_1 =
\<f, (-L)f\>_\rho = \<f, (-S)f\>_\rho$.  The $H_1$ space then is the
completion with respect to this norm.  With respect to
``coefficient'' operators, we have the corresponding $H_1$
(semi)-norm on functions $\f$ supported on $\Epsilon$ 
given by $\|\f\|_1^2 = \<\f,(-\L)f\> = \|f\|_1^2$, and corresponding
completed space $H_1$. 

Let $H_{-1}$ be the dual of $H_1$, namely, the completion over local functions with respect to norms
$\|\cdot\|_{-1}$ given by
\begin{eqnarray*}
\|f\|^2_{-1} & = & \sup_g \big\{2\<f,g\>_\rho -
\|g\|_1^2\big\}\\
&=&\sup_{\g} \big\{2\<\f,\g\> -
\|\g\|_1^2\big\} \ \ = \ \|\f\|^2_{-1}.
\end{eqnarray*}    

Similarly, we define, for convenience, the notation
$\|f\|^2_{1,\lambda} = \|\f\|^2_{1,\lambda}=\<f,(\lambda-S)f\>_\rho$ and
$\|f\|^2_{-1,\lambda} = \|\f\|^2_{-1,\lambda}=\sup_g \{2\<f,g\>_\rho - \|g\|^2_{1,\lambda}\}$.

Then, in this notation, we write for local $f$ that
\begin{eqnarray*}
\<f,(\lambda - L)^{-1}f\>_\rho &=& \sup_g \big\{2\<f,g\>_\rho
-\|g\|_{1,\lambda}^2 -\|Ag\|_{-1,\lambda}^2\big\}\\
&=& \sup_\g \big\{2\<\f,\g\>
-\|\g\|_{1,\lambda}^2 -\|\A\g\|_{-1,\lambda}^2\big\}.
\end{eqnarray*}

\vskip .2cm

{\it ``Free Particle'' Bounds.}
To analyze these variational formulas, it will be helpful
computationally to ``remove the hard-core exclusion.''  In other
words, we want to get equivalent bounds in terms of operators which
govern completely independent or ``free'' motions.  We follow
Bernardin \cite{B}.  Let $\chi_n = (Z^2)^n$ and note that $\Epsilon_n \subset
\chi_n$.
Consider $n$ independent random walks with jump rates $s$ on
$Z^2$.  The process $x_t = (x^1_t,\ldots,x^n_t)$ evolves on $\chi_n$
and has generator acting on finitely supported functions
$$(\S_{{\rm free}}\f)(x) \ = \ \sum_{1\leq j\leq n \atop z\in
  Z^2}s(z)[\f(x+ze_j)-\f(x)].$$
With respect to finitely supported functions on $\chi_n$, let 
$$\<\phi,\psi\>_{{\rm free}} = \frac{1}{n!}\sum_{x\in \chi_n}
\phi(x)\psi(x)$$
be the innerproduct, and denote the norms $\|\f\|_{1,{\rm free}}$ and
$\|\f\|_{-1,{\rm free}}$
by
$$\begin{array}{rl}
\|\f\|^2_{1,{\rm free}} & = \ \frac{1}{n!}\sum_{x\in
  \chi_n}\sum_{1\leq j\leq n \atop z\in Z^d}[\f(x+ze_j)-\f(x)]^2\\
\|\f\|^2_{-1,{\rm free}}&= \ \sup_{\g}\big\{2\<\f,\g\>_{{\rm free}} -
\|\g\|^2_{1,{\rm free}}\big\}.
\end{array}
$$
Define also the resolvent quantities
$$\begin{array}{rl}
\|\f\|^2_{1,\lambda, {\rm free}} & = \ \<\f,(\lambda-\S_{{\rm
    free}})\f\>_{{\rm free}}\\
\|\f\|^2_{-1,\lambda, {\rm free}}&= \ \sup_{\g}\big\{2\<\f,\g\>_{{\rm free}} -
\|\g\|^2_{1,\lambda, {\rm free}}\big\}.
\end{array}
$$
 
Let $\G_n\subset \chi_n$ be those points whose coordinates are
distinct.  The following is a part of Theorems 3.1 and 3.2 \cite{B}
[which simplifies as $\tilde{\f} = f$ for $f\in \C_1$, and $1_{x\in \G_n}
\tilde{\f} = \f$ for $\f\in \C_n$].
\begin{prop}
\label{free}
There exists a constant $C_5=C_5(n)$ such that for all functions
in $\C_1$ we have
$$C_5^{-1}\|\f\|^2_{1,{\rm free}} \ \leq \ \|\f\|_1^2 \ \leq C_5\|\f\|^2_{1,{\rm free}}.$$
Also, for all functions in $\C_n$ (for any $n$),
$$ \|\f\|^2_{-1,\lambda} \ \leq \ C_5\| \f\|^2_{-1, \lambda,
  {\rm free}}.$$
\end{prop}

We express now the ``free'' $H_1$ and $H_{-1}$ norms in terms of
Fourier transforms.  Let $\psi$ be a local function on $\chi_n$ and let $\widehat{\psi}$ be
its Fourier transform
$$\widehat{\psi}(s_1,\ldots,s_n) = \frac{1}{\sqrt{n!}} \sum_{x\in \chi_n}
e^{2\pi i (x_1s_1 + \cdots x_ns_n)}\psi(x)$$
where $s_1,\ldots,s_n\in [0,1]^2$.
Compute $\widehat{\S_{\rm free}}$ from the relation $\widehat{\S_{\rm free}}\widehat{\psi} =
\widehat{\S_{{\rm free}}\psi}$ as
$$\widehat{\S_{\rm free}}\widehat{\psi}(s_1,\ldots,s_n) \ = \ -\bigg[\sum_{j=1}^n
\theta_2(s_j)\bigg] \widehat{\psi}(s_1,\ldots,s_n)$$
where $\theta_2(u) = 2\sum_{z\in Z^2} s(z)\sin^2(\pi(u\cdot z))$.

Hence,
$$\|\psi\|^2_{1,\lambda, {\rm free}} \ = \ \int_{s\in ([0,1]^2)^n}
\bigg(\lambda + \sum_{j=1}^n
\theta_2(s_j)\bigg)|\widehat{\psi}(s_1,\ldots,s_n)|^2 ds$$
and
$$\|\psi\|^2_{-1,\lambda, {\rm free}} \ = \ \int_{s\in ([0,1]^2)^n}
\frac{|\widehat{\psi}(s_1,\ldots,s_n)|^2}{\lambda + \sum_{j=1}^n
\theta_2(s_j)} ds.$$

\subsect{Proof of Proposition \ref{lb1}}
Let $f(\eta) = \eta_0 - 1/2$.  Then, $\f =
(1/2)\delta_0$ where $\delta_0$ is the indicator of the
set $\{0\}$.  To prove Proposition \ref{lb1}, we find lower bounds on
$\|\delta_0\|^2_{-1,\lambda}$.  From Proposition \ref{comp}, we will
assume $L$ takes form $L=L_0$.
Write now, as $\rho = 1/2$, that
$$
\<\delta_0,(\lambda - \L)^{-1}\delta_0\> \ = \ \sup_{\phi} \bigg\{2\<\delta_0,\phi\> - \|\phi\|_{1,\lambda}^2 -
  \|\A\phi\|_{-1,\lambda}^2\bigg\}.$$
The strategy now will be (1) to replace to restrict the supremum on
$\phi$ to local
degree $1$ functions in $\C_1$,
and (2) 
to use the comparison bounds with respect to
independent particles (Proposition \ref{free}) to help bound terms in the
variational formula.

Let now $\phi$ be a degree one function.
To simplify notation, let
$s(\pm e_1) = b_1>0$, $s(\pm e_2)=b_2>0$ and $a(e_1) =
-a(-e_1)=a_1\neq 0$, $a(e_2)=-a(-e_2)=a_2$.
Note now that
$\A\phi$ takes form $\A\phi = \A_{1,2}\phi$.  More specifically,
$(\A_{1,2}\phi)(x,y) = a(y-x)(\phi(x)-\phi(y))$ is supported on
distinct two-tuples and can be written as
$$(\A_{1,2}\phi)(x,y) \ =\ \left\{\begin{array}{rl}
a_1(\phi(x) - \phi(x+e_1)) & \ {\rm when \ }y = x+e_1\\
a_2(\phi(x) - \phi(x+e_2)) & \ {\rm when \ }y = x+e_2\\
0&\ {\rm otherwise.}\end{array}\right.
$$

Using Proposition \ref{free}, we have for some constant $C_6$ that
\begin{equation}
\label{vari}
\<\delta_0,(\lambda - \L)^{-1}\delta_0\> \ \geq \ 
C_6\sup_{\phi \ {\rm local \ in \ } \C_1} \bigg\{2\<\delta_0,\phi\> -
\|\phi\|_{1,\lambda, {\rm free}}^2 -
  \|\A_{1,2}\phi\|_{-1,\lambda, {\rm free}}^2\bigg\}.
\end{equation}

Now, it is a calculation to find that
\begin{eqnarray*}
\widehat{\A_{1,2}\phi}(s,t) & = & \frac{1}{\sqrt{2}}\sum_{(x,y)\in
  (Z^d)^2} e^{2\pi i (x\cdot s + y\cdot t)}(\A_{1,2}\phi) (x,y)  \\
&=& \frac{i}{\sqrt{2}}\widehat{\phi}(s+t)[2a_1\sin(2\pi s_1)
  +2a_2\sin(2\pi s_2)\\
&&\ \ \ \ \ \ \ \ \ \ \ \ \ \ \ \ \ \ \ \ \ \ \ \ \ \ \ \ \ 
  +2a_1\sin(2\pi t_1) + 2a_2\sin(2\pi t_2)].
\end{eqnarray*}
Then, the expression in brackets in (\ref{vari}) becomes in Fourier terms
\begin{eqnarray*}
&&\int_{[0,1]^2}\bigg(\widehat{\phi}(s) - (\lambda +
\theta_2(s))|\widehat{\phi}|^2\bigg) ds\\
&&-\frac{1}{2}\int_{([0,1]^2)^2} \frac{\big (\sum_{i=1}^2 2a_i\sin(2\pi s_i) +2a_i\sin(2\pi t_i)
  \big )^2}{\lambda +\theta_2(s)
  +\theta_2(t)} |\widehat{\phi}(s+t)|^2 ds_1ds_2 dt_1dt_2.
\end{eqnarray*}

We now change coordinates in the second integral:
$$(s_1,s_2,t_1,t_2) \ = \ (\frac{u+v}{2}, \frac{w+z}{2},
\frac{u-v}{2}, \frac{w-z}{2})$$
(whose Jacobian determinant in absolute value is $1/4$).  The region
$[0,1]^4$ is mapped to $D^2$ where $D$ is a planar diamond with
vertices $(0,0), (1,-1), (1,1), (2,0)$.
Let
\begin{eqnarray}
\label{gamma}
\gamma(u,v,w,z) &=& 2b_1\sin^2(\pi\frac{u+v}{2}) + 2b_2\sin^2(\pi
\frac{w+z}{2}) \\
&&\ \ \ \ \ \ +2b_1\sin^2(\pi \frac{u-v}{2}) + 2b_2\sin^2(\pi \frac{w-z}{2})\nonumber\\
&=&
4b_1\sin^2(\pi(u/2))\cos^2(\pi(v/2)) +
4b_1\sin^2(\pi(v/2))\cos^2(\pi(u/2))\nonumber \\
&&\ \ \ \ \ \ 
+4b_2\sin^2(\pi(w/2))\cos^2(\pi(z/2)) +4b_2\sin^2(\pi(z/2))\cos^2(\pi(w/2))\nonumber
\end{eqnarray}
and
\begin{eqnarray*}
\upsilon(u,v,w,z) & = & 4a_1^2\sin^2(\pi u)\cos^2(\pi v) +4a_2^2\sin^2(\pi
  w)\cos^2(\pi z)\\
&&\ \ \ 
  +8a_1a_2\sin(\pi u)\cos(\pi v)\sin(\pi w)\cos(\pi z)\\
&\leq& 8a_1^2\sin^2(\pi u)\cos^2(\pi v) +8a_2^2\sin^2(\pi
  w)\cos^2(\pi z).
\end{eqnarray*}
The second integral is 
rewritten as
\begin{eqnarray*}
&&\int_D\int_{D} \frac{\upsilon(u,v,w,z)}{\lambda + \gamma(u,v,w,z)}
|\widehat{\phi}(u,w)|^2 dudvdwdz\\
&&\ \ \ \ \ \ \leq \ \int_D\int_{D} \frac{8a_1^2\sin^2(\pi u)\cos^2(\pi v) +8a_2^2\sin^2(\pi
  w)\cos^2(\pi z)}{\lambda + \gamma(u,v,w,z)}
|\widehat{\phi}(u,w)|^2 dudvdwdz .
\end{eqnarray*}

By changing variables and adding and identifying some parts of the
region of integration, it is not difficult to see that
the last integral reduces to
\begin{eqnarray*}
&& 4\int_0^{1}\int_0^{1}\int_0^1\int_0^1 
\frac{8a_1^2\sin^2(\pi u)\cos^2(\pi v) +8a_2^2\sin^2(\pi
  w)\cos^2(\pi z)}{\lambda + \gamma(u,v,w,z)}
|\widehat{\phi}(u,w)|^2 dvdzdudw\\
&&\ \ = 4\int_0^1\int_0^1 [8a_1^2\sin^2(\pi u)F^1_\lambda(u,w) + 8a_2^2\sin^2(\pi w)F^2_\lambda (u,w)]|\widehat{\phi}(u,w)|^2 dudw
\end{eqnarray*}
where
$$F^1_\lambda (u,w) = \int_0^1\int_0^1 \frac{\cos^2(\pi v)}{\lambda + \gamma(u,v,w,z)}
dv dz \ \ {\rm and \ \ } F^2_\lambda (u,w) = \int_0^1\int_0^1 \frac{\cos^2(\pi z)}{\lambda + \gamma(u,v,w,z)}
dv dz.$$

Let now
$C_7 =
2(b_1^2 +b_2^2) +16(a_1^2 +a_2^2)$.
For general rates, substituting into (\ref{vari}), we obtain
$C_6^{-1}\<\delta_0, (\lambda -\L)^{-1}\delta_0\>$ greater than
\begin{eqnarray*}
&&\sup_{\phi}\bigg\{ \int_0^1\int_0^1 
\widehat{\phi}(u,w)  \\
&&\ \ \ - \big(\lambda + C_7(\sin^2(\pi u) +\sin^2(\pi w))(1+F^1_\lambda (u,w)
+F^2_\lambda(u,w))
\big ) |\widehat{\phi}(u,w)|^2 du dw \bigg \}
\end{eqnarray*}
where the supremum is on $\phi$ local, or without loss of generality
on $L^2(\Epsilon_1)$.  
When $a_2=0$, we have the lower bound
\begin{eqnarray}
\label{special}
&&\sup_{\phi}\bigg\{ \int_0^1\int_0^1 
\widehat{\phi}(u,w)   \\
&&\ - \big(\lambda + C_7(\sin^2(\pi u) +\sin^2(\pi w))
+C_7\sin^2(\pi u)F^1_\lambda (u,w)
\big ) |\widehat{\phi}(u,w)|^2 du dw \bigg \}.\nonumber
\end{eqnarray}

We now concentrate on the general rates case.
By optimizing on $\phi$ we get the lower bound
\begin{equation}
\label{optimized}
\frac{1}{4}\int_0^1\int_0^1 \frac{du dw}{\lambda + C_7(\sin^2(\pi u)
  +\sin^2(\pi w))(1+F^1_\lambda (u,w)+F^2_\lambda(u,w))}
\end{equation}
with optimizer
$$\widehat{\phi}(u,w) \ = \ \frac{1}{2} \frac{1}{\lambda +
  C_7(\sin^2(\pi u) +\sin^2( \pi w))(1+F^1_\lambda (u,w)+F^2_\lambda(u,w)) }$$
which is the transform of a real function as $\widehat{\phi}(u,w) =
  \widehat{\phi}^*(1-u,1-w)$ (note
  $F_\lambda^i(u,w)=F^i_\lambda(1-u,1-w)$ for $i=1,2$ by observing
  (\ref{gamma}) and changing variables $v\rightarrow 1-v$ and
  $z\rightarrow 1-z$).


We now bound $F^1_\lambda(u,w) +F^2_\lambda(u,w)$ for $|u|,|w|\leq 1/2$.
Since $\cos(x)$ is decreasing for $0\leq x\leq \pi/2$ and 
$\sin(x) \geq (2/\pi)x$ for $0\leq x\leq \pi/2$, 
we have
\begin{eqnarray*}
F^1_\lambda (u,w) +F^2_\lambda(u,w)
&\leq& \int_{[([1/2,1]\times [0,1]) \cup ([0,1]\times [1/2,1])} \frac{2dv dz}{\lambda
+2b_1v^2 +2b_2z^2}\\
&&\ \ \ \ \ + \int_0^{1/2}\int_0^{1/2} \frac{2dv dz}{\lambda
+2b_1(u^2 +v^2) +2b_2(w^2+ z^2)}\\
&\leq &C_8 + \int_0^{\pi/2}\int_0^{1} \frac{2s ds
  d\alpha}{\lambda + \b(u^2 +w^2 +s^2)}\\
&=& C_8+\frac{\pi}{2\b} \log\big[\frac{\lambda + \b(u^2 +w^2) +\b}{\lambda + \b(u^2
  +w^2)}\big ]\\
& \leq & C_9 + \frac{\pi}{2\b}|\log (\lambda +  \b(u^2 +w^2))|.
\end{eqnarray*}
where $\b = 2\min\{b_1,b_2\}$ and $C_8=C_8(b_1,b_2)$, $C_9=C_9(b_1,b_2)$ are constants.

Hence, as $\sin(x) \leq x$, we can bound (\ref{optimized}) below by
\begin{equation}
\label{afterF}
\frac{1}{4}\int_0^{1/2}\int_0^{1/2} \frac{du dw }{\lambda +
  4C_7\pi^2 (u^2+w^2)(1 + C_9 + \frac{\pi}{2\b}|\log(\lambda + \b
  (u^2+w^2))|)}.
\end{equation}
We have with respect to constants $C_{10}, C_{11}$ that
$4$ times (\ref{afterF}) is greater than
\begin{eqnarray*}
\int_0^{1/2} \frac{r dr }{\lambda +
  C_{10}r^2(1 + |\log(\lambda + C_{10}r^2)|)}
&=& \int_0^{1/(2\sqrt{\lambda})} \frac{rdr }{1+
  C_{10}r^2(1+ |\log\lambda(1+ C_{10}r^2)|)}\\
&\geq& \frac{1}{C_{11}}\int_1^{1/(2\sqrt{\lambda})} \frac{ dr}{
r|\log\lambda r^2|}\\
&\geq& \frac{1}{C_{11}}\int_{\sqrt{\lambda}}^{1/2}
  \frac{dr}{r|\log r^2|}.
\end{eqnarray*}
This last expression is order $|\log(\log \lambda)|$.  To get the
larger expected
order of $|\log \lambda|^{2/3}$, it seems one would need to optimize also over higher degree
functions in (\ref{vari}). 

We note in the case $a_2=0$, we bound (\ref{special}) by
$$
\frac{1}{4}\int_0^{1/2}\int_0^{1/2} \frac{du dw }{\lambda +
  4C'_7\pi^2 (u^2+w^2) + C'_7u^2(C_9 + \frac{\pi}{2\b}|\log(\lambda + \b
  (u^2+w^2))|)}.
$$
Following closely the sequence to bound the second-class
particle variance in $d=2$ (cf. p. 470 \cite{LQSY}),
we observe $|\log(\lambda +u^2+w^2)| \leq |\log(\lambda +w^2)|$ 
for $\lambda$
small and $0\leq
u,w\leq 1/2$.  And so, we obtain a lower bound on order
$$
\int_0^{1/2}\int_0^{1/2} \frac{du dw}
{\lambda +u^2+w^2 +u^2|\log (\lambda + w^2)|}.$$
With substitution $u = y(1+|\log(\lambda +w^2)|)$ the above expression
is bounded below by
$$
\int_0^{1/2}\int_0^{1/2} \frac{dy dw}
{\lambda +y^2+w^2}(1+|\log (\lambda + w^2)|)^{-1/2}.$$
Changing to polar coordinates and restricting $\pi/6\leq \alpha \leq
\pi/4$,
we get a lower bound on order as in \cite{LQSY}
$$\int_0^{1/20}\frac{rdr}{\lambda +r^2} |\log(\lambda +r^2)|^{-1/2} \
\geq \ C_{12} |\log \lambda |^{1/2}$$
for a constant $C_{12}$.
\qed

\vskip .5cm
{\bf Acknowledgements.}  I would like to thank C. Bernardin,
C. Landim, S. Olla, J. Quastel, H. Spohn, and H.T. Yau for
useful conversations.
\vskip .2cm

\bibliographystyle{plain}

\vskip.5cm
\noindent Sunder Sethuraman\\
400 Carver Hall\\
Dept. of Mathematics\\
Iowa State University\\
Ames, IA \ 50011\\
sethuram@iastate.edu

\end{document}